\documentclass[11pt]{article}
\usepackage{amsmath,amsfonts,here,latexsym,amssymb,annals3,eufrak,amscd, graphicx,pictexwd,pstricks,pst-plot}
\usepackage[square,numbers]{natbib}
\def\a{\alpha}
\def\b{\beta}
\def\R{\mathbb R}

\def\labda1{\lambda_1}
\def\labda2{\lambda_2}

\def\f{\phi}

\def\s{\sigma}

\def\comment#1{\relax}

\def\=in{\mathop{\rm =}}

\def\eop{\hfill\mbox{$\Box$}\newline}

\newtheorem{theorem}{Theorem}

\newtheorem{definition}{Definition}

\renewcommand{\a}{\alpha}

\begin{document}
\title{The remaining area of the convex hull of a Poisson process}
\author{Piet Groeneboom}
\date{\today}
\affiliation{Delft University of Technology}
\AMSsubject{5A22,52A38}
\maketitle
\begin{abstract}
In \cite{cabogr:94} the remaining area of the left-lower convex hull of a Poisson point process with intensity one in the first quadrant of the plane was analyzed, using the methods of \cite{gr:88}, giving formulas for the expectation and variance of the remaining area for a finite interval of slopes of the boundary of the convex hull. However, the time inversion argument of 
\cite{gr:88} was not correctly applied in \cite{cabogr:94}, leading to an incorrect scaling constant for the variance.
The purpose of this note is to show how the correct application of the time inversion argument gives the right expression, which is in accordance with results in \cite{nagaev_kham:91} and \cite{buchta:03}.
\end{abstract}

\section{Introduction}
\label{sec:intro}
Let ${\cal P}$ be a homogeneous Poisson point process of intensity $1$ in the first quadrant of the plane. We consider the left-lower convex hull of this point process, as in \cite{gr:88}, where a Markovian jump process of points $W(a)$ is introduced, by which one can ``walk" through the vertices of the convex hull. The points $W(a)$ are defined in the following way.

\begin{definition}
\label{def:W(A)}
{\rm
For each $a>0$, $W(a)=(U(a),V(a))$ is the point of the realization of the Poisson process ${\cal P}$ on $\R_+^2$ such that all points of the
realization lie to the right of the line of the line $x+ay=c$ which passes through $W(a)$. If there are several of such points (which happens with
probability zero for fixed $a$), we take the point with the smallest $y$-coordinate.
}
\end{definition}

In \cite{piet:11a} the following result is proved for the joint distribution of the number of jumps and the remaining area, corresponding to a (slopes) interval $[a,b]$.

\begin{theorem}
\label{th:CL_for_N_D}
Let $N(a,b)$ be the number of jumps in the interval $[a,b]$ of the process $W$, as defined in Definition \ref{def:W(A)}, and let $D(a,b)$ be the area of the union of the triangles $T_i$, corresponding to points of jump $a_i\in[a,b]$, as defined in Theorem
2.1 of \cite{piet:11a}. Then:
$$
\left(\tfrac5{27}\log(b/a)\right)^{-1/2}\left(N(a,b)-\tfrac13\log(b/a),D(a,b)-\tfrac13\log(b/a)\right)\stackrel{{\cal D}}\longrightarrow
N(0,\Sigma),\,b/a\to\infty,
$$
where $N(0,\Sigma)$ is a bivariate normal distribution with expectation $0$ and covariance matrix defined by
\begin{equation}
\label{def_sigma}
\Sigma=\left(\begin{array}{ll}
1 &1\\
1 &\displaystyle{\tfrac{14}{5}}
\end{array}
\right)
\end{equation}
\end{theorem}

This result is used in \cite{piet:11a} to prove the unpublished result in \cite{nagaev_kham:91}, which gives a central limit theorem for the joint distribution of the number of vertices and the remaining area of the convex hull of a uniform sample of points from the interior of a convex polygon, see Theorem 1.1 in \cite{piet:11a}. This extends the results in \cite{cabogr:94} and \cite{gr:88}.

Theorem \ref{th:CL_for_N_D} is illustrated in \cite{piet:11a} by a simulation of the jump process, generating the vertices and the areas, without first having to generate the points of the Poisson process for which one considers the convex hull. In this simulation Theorem \ref{theorem_Poisson} below is used. It is the purpose of this note to provide a proof of this result, which at the same time corrects the proof in \cite{cabogr:94}, which went astray by an incorrect application of a time inversion argument, as explained in section \ref{sec:area}.

\section{Expectation and variance of the area}
\label{sec:area}

In this section, we prove the following theorem.

\begin{theorem}
\label{theorem_Poisson}
Let, for $0<a<b<\infty$, $D(a,b)$ be the area of the region,
bounded on the right and left by vertical lines through the $x$-coordinates of the points $W(a)$ and $W(b)$, and bounded from below and
above by the line $y=0$ and the (left lower) boundary of the convex hull of $\cal P$, respectively. Morever, let, for $0<a<b<\infty$,
$\b=b/a$ and $\a=\b-1$. Then
\begin{enumerate}
\item[(i)] $E D(a,b)=\tfrac13\log\b$,
\item[(ii)]
\begin{align*}
&\mbox{\rm Var}(D(a,b))\\
&=\frac{14}{27}\log\b+\frac2{3\a^2}+\frac4{9\a}-\frac{44}{45}
-\frac{2\{3+\a(3-4\a)\}\tan^{-1}\left(\sqrt{\a}\right)}{9\a^{5/2}}
+\frac49\left(\tan^{-1}\left(\sqrt{\a}\right)\right)^2.
\end{align*}
\end{enumerate}
\end{theorem}

\noindent
{\bf Proof.} We follow the line of the argument in Appendix A of \cite{cabogr:94}. First of all, we get, as in (4.3)
of \cite{cabogr:94}, replacing $A(a,b)$ by $D(a,b)$,
\begin{align*}
\lim_{h\downarrow0}h^{-1}E\left\{D(a,a+h)|W(a)\right\}=\tfrac5{24}aV(a)^4,
\end{align*}
and hence
$$
ED(1,a)=\tfrac5{24}\int_1^a sEV(s)^4\,ds=\tfrac5{24}\int_1^a s^{-1}EV(1)^4\,ds=\tfrac13\log a.
$$
Stationarity in the logarithmic scale then also gives $ED(a,b)=\tfrac13\log(b/a)$ and (i) follows.

To prove (ii), we start with relation (4.4) in \cite{cabogr:94}, which is repeated below (with $A(1,a)$ replaced by $D(1,a)$):
\begin{equation}
\label{second_moment}
ED(1,a)^2=\tfrac{11}{120}\int_1^a s^2EV(s)^6 ds+\tfrac5{12}\int_1^a E\left\{D(1,s)V(s)^4\right\}\,ds.
\end{equation}
As noted on the top of page 46 of
\cite{cabogr:94}, the first integral is equal to $(22/35)\log a$.
As indicated on page 46 of \cite{cabogr:94}, the other integral can be computed by a time reversal argument. But this is also exactly
the point where things went astray.

\begin{figure}
\psset{unit=5mm} 
\begin{center}
\begin{pspicture}(0,0)(10,5) 
\psaxes[labels=none,ticks=none,showorigin=true](0,0)(10,5)
\uput[d](-1,0.5){(0,0)}
\pspolygon[fillstyle=solid,fillcolor=lightgray](2,0)(4,0)(4,2.5)(2,5)
\pspolygon[fillstyle=solid,fillcolor=lightgray](4,0)(7,0)(7,1)(4,2.5)
\pspolygon[fillstyle=solid,fillcolor=lightgray](7,0)(10,0)(10,0.5)(7,1)
\uput[d](6,6.8){$W(1)=(U(1),V(1))$}
\uput[d](14,2.2){$W(s)=(U(s),V(s))$}
\uput[d](6.2,3.8){$D(1,s)$}
\psline[arrows=->](2.5,5.8)(2,5)
\psline[arrows=->](6,2.5)(5,1)
\psline[arrows=->](10.5,1.2)(10,0.5)
\end{pspicture}
\caption{$D(1,s)=$ area of shaded region}
\label{figure1}
\end{center}
\end{figure}
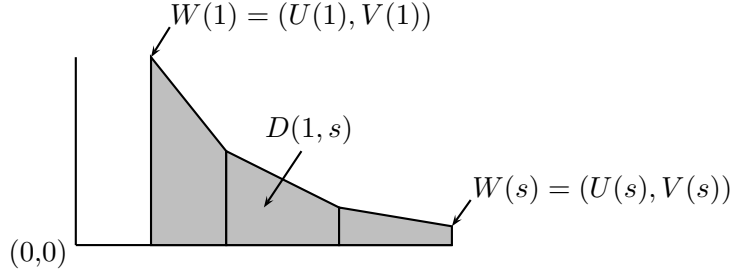

\begin{figure}
   \psset{unit=5mm} 
   \begin{center}
      \begin{pspicture}(0,-0.5)(5,12) 
						\psaxes[labels=none,ticks=none,showorigin=true](0,0)(5,10)
						\uput[d](-1,0.5){(0,0)}
  				\pspolygon[fillstyle=solid,fillcolor=lightgray](0,2)(0,4)(2.5,4)(5,2)
						\pspolygon[fillstyle=solid,fillcolor=lightgray](0,4)(0,7)(1,7)(2.5,4)
					\pspolygon[fillstyle=solid,fillcolor=lightgray](0,7)(0,10)(0.5,10)(1,7)
					\uput[d](6.2,11.0){$W(1/s)=(U(1/s),V(1/s))$}
					\uput[d](9,4){$W(1)=(U(1),V(1))$}
\uput[d](4.3,7){$D'(1/s,1)$}
\psline[arrows=->](2.5,6)(1,5)
\psline[arrows=->](1.5,10.2)(0.5,10)
\psline[arrows=->](5.5,2.8)(5,2)
      \end{pspicture}
      \caption{$D'(1/s,1) =$ area of shaded region}
      \label{figure2}
   \end{center}
\end{figure}
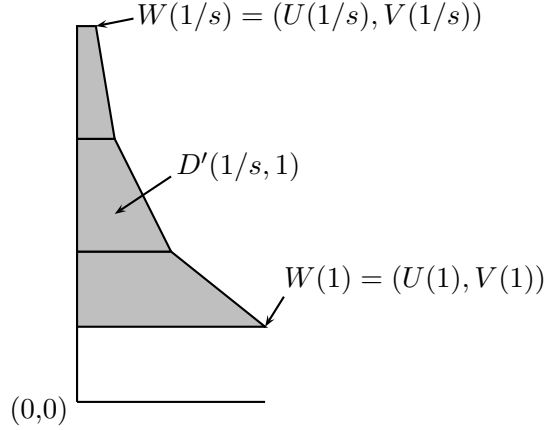

We have to compute the expectation of the product $D(1,s)V(s)^4$. $D(1,s)$ is area of a concatenation of 
trapezia, as for example shown in Figure \ref{figure1}. Note that Figure \ref{figure1} corresponds to three successive values of the jump
process $s\mapsto D(1,s)$. For reasons of symmetry it is clear that $E D(1,s)V(s)^4$ will be equal to
$E D'(1/s,1)U(1/s)^4$, where $D'(a,b)$ is the area of the region, bounded from above and below by horizontal lines through the
$y$-coordinates of the points $W(a)$ and $W(b)$, and bounded on the left and right by the $y$-axis and the (left lower) boundary of the
convex hull of $\cal P$, respectively, see Figure \ref{figure2}. The advantage of working with
$D'(1/s,1)U(1/s)^4$ instead of $ D(1,s)V(s)^4$ is that the somewhat elusive concatenation of trapezia, corresponding
to the ``developing area" (to use a phrase, coined in \cite{cabogr:94}), now lies {\it in the future} of the process with respect to
$U(1/s)$, whereas otherwise, in working directly with $D(1,s)V(s)^4$, we would have to condition on the future value
$V(s)$ of the process.

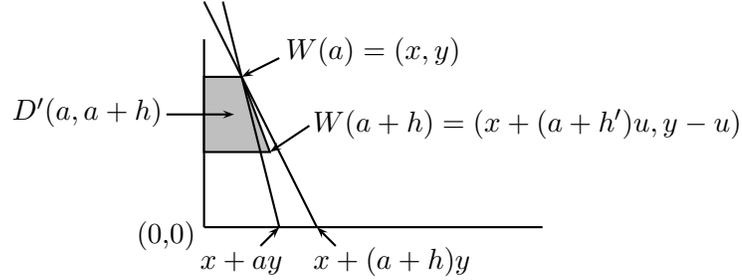
\begin{figure}
\psset{unit=5mm} 
\begin{center}
\begin{pspicture}(-3,-1)(8,5) 
\psaxes[labels=none,ticks=none,showorigin=true](0,0)(9,5)
\uput[d](-1,0.5){(0,0)}
\pspolygon[fillstyle=solid,fillcolor=lightgray](0,4)(1,4)(1.75,2)(0,2)
\psline[](0,6)(3,0)
\psline[](0.5,6)(2,0)
\uput[d](4.5
,5.4){$W(a)=(x,y)$}
\uput[d](8.6,3.5){$W(a+h)=(x+(a+h')u,y-u)$}
\uput[d](1,-0.3){$x+ay$}
\uput[d](5,-0.18){$x+(a+h)y$}
\uput[d](-3.1,3.8){$D'(a,a+h)$}
\psline[arrows=->](2,4.5)(1,4)
\psline[arrows=->](2.8,2.7)(1.75,2)
\psline[arrows=->](1.5,-0.5)(2,0)
\psline[arrows=->](3.5,-0.5)(3,0)
\psline[arrows=->](-1,3)(0.8,3)
\end{pspicture}
\caption{The increment $D'(a,a+h)$ in the time interval $[a,a+h]$, Note: $h'\in[0,h]$.}
\label{figure3}
\end{center}
\end{figure}

We have:
\begin{align*}
&\lim_{h\downarrow0}h^{-1}E\left\{D'(a,a+h)|W(a)=(x,y)\right\}\\
&=\int_0^y u\left(ux+\tfrac12au^2\right)du\\
&=\tfrac13xy^3+\tfrac18ay^4.
\end{align*}
The situation is illustrated in Figure \ref{figure3}, where, during the time interval $[a,a+h]$, a new point $W(a+h)$ is caught, leading to
the increment $D'(a,a+h)$, which is a trapezium, bounded on the left by the $y$-axis. If $W(a)=(x,y)$, the new point $W(a+h)$ has
coordinates $(x+(a+h')u,y-u)$, where $h'\in[0,h]$ and $u\in(0,y)$.

Hence, using the time inversed process $a\mapsto R(a)=(S(a),T(a))=(V(1/a),U(1/a))$, as on page 343 of
\cite{gr:88}, we get:
\begin{align*}
E\{D(1,a)|W(a)\}&=E\{D'(1/a,1)|W(1/a)\}\\
&=\int_{1/a}^1 \left\{\tfrac13E\left\{S(s)T(s)^3\bigm|R(1/a)\right\}+\tfrac18E\left\{sT(s)^4\bigm|R(1/a)\right\}\right\}\,ds\\
&=\int_1^a s^{-2}\left\{\tfrac13E\left\{S(1/s)T(1/s)^3\bigm|R(1/a)\right\}+\tfrac18E\left\{s^{-1}T(1/s)^4\bigm|R(1/a)\right\}\right\}\,ds\\
&=\int_1^a s^{-3}\left\{\tfrac13 sE\left\{V(s)U(s)^3\bigm|W(a)\right\}+\tfrac18E\left\{U(s)^4\bigm|W(a)\right\}\right\}\,ds,
\end{align*}
and, instead of $(5/12)EU(s)^4V(a)^4$, we have to compute:
\begin{equation}
\label{second_integrand}
\tfrac13 s EU(s)^3V(s)V(a)^4+\tfrac18EU(s)^4V(a)^4.
\end{equation}
Note that the mistake in computing the second integral in (4.4) of \cite{cabogr:94} was caused by not taking into account that the relevant
trapezium with upper right vertex $W(1/s)$ is bounded on the left side by the
$y$-axis and bounded below by a horizontal line above the $x$-axis, as is shown (for several of such trapezia) in Figure \ref{figure2},
instead of being bounded below by the $x$-axis and bounded on the left by a vertical line to the right of the $y$-axis as in Figure
\ref{figure1}. 

As in \cite{cabogr:94}, we can use the stationarity in the logarithmic scale to reduce the computation of (\ref{second_integrand}) to
the computation of
\begin{equation}
\label{second_integrand2}
\tfrac13EU(1)^3V(1)V(s)^4+\tfrac18EU(1)^4V(s)^4,\,s\ge1.
\end{equation}
This expectation is computed in the the appendix. The result is:
\begin{align}
\label{relevant_exp}
&\tfrac13EU(1)^3V(1)V(s)^4+\tfrac18EU(1)^4V(s)^4
=\frac{210+2\s\{85+4(\s-1)\s\}}{15\s^4(1+\s)}-\frac{2(7+\s)\tan^{-1}\left(\sqrt{\s}\right)}{\s^{9/2}}\,,
\end{align}
where $\s=s-1$. We note in passing that the right side of (\ref{relevant_exp}) tends to $104/315$, as $\s\downarrow0$, so there is no singularity at
$\s=0$, as is also clear from the expression at the left side of (\ref{relevant_exp}). Let
$$
\f(s)=\tfrac13EU(1)^3V(1)V(s)^4+\tfrac18EU(1)^4V(s)^4,\,s\ge1.
$$
Then we can write, using (\ref{relevant_exp}),
$$
\f(s)=\frac{24+2s\{105+4s(s-4)\}}{15s(s-1)^4}-\frac{2(6+s)\tan^{-1}\left(\sqrt{s-1}\right)}{(s-1)^{9/2}}\,,\,s\ge1.
$$

Thus the computation of the second integral in (\ref{second_moment}) boils down to the computation of
$$
\tfrac5{12}\int_1^a s\,ds\int_1^s\frac{\f(s/r)}{r^3}\,dr=\tfrac5{12}\int_1^a \frac{ds}{s}\int_1^s u\f(u)\,du.
$$
We have:
\begin{align*}
&\tfrac5{12}\int_1^s u\f(u)\,du
=\frac{4-19s+2(s-1)^3\log s}{9\s^3}+\frac{\{s(14+5s)-4\}\tan^{-1}\left(\sqrt{\s}\right)}{9\s^{7/2}}-\frac{104}{945}\,.
\end{align*}
Thus
\begin{align*}
&\tfrac5{12}\int_1^a \frac{ds}{s}\int_1^s u\f(u)\,du\\
&=\tfrac19(\log a)^2-\tfrac{104}{945}\log a+\tfrac49\left(\tan^{-1}\left(\sqrt{a-1}\right)\right)^2\\
&\qquad\qquad\qquad\qquad\qquad\qquad+\frac{2\{4+a(4a-11)\}\tan^{-1}\left(\sqrt{a-1}\right)}{9(a-1)^{5/2}}+\frac{2+4a}{9(a-1)^2}-\frac{44}{45}\,.
\end{align*}
Using the notation $\a=a-1$, this can also be written
\begin{align*}
&\tfrac5{12}\int_1^a \frac{ds}{s}\int_1^s u\f(u)\,du\\
&=\tfrac19(\log a)^2-\tfrac{104}{945}\log a+\tfrac49\left(\tan^{-1}\left(\sqrt{\a}\right)\right)^2
-\frac{2\{3+\a(3-4\a)\}\tan^{-1}\left(\sqrt{\a}\right)}{9\a^{5/2}}+\frac{6+4\a}{9\a^2}-\frac{44}{45}\,.
\end{align*}
Hence, taking into account that the first term on the right side of (\ref{second_moment}) equals $(22/35)\log a$, we get that the right side of
(\ref{second_moment}) is equal to:
\begin{align*}
\tfrac19(\log a)^2+\tfrac{14}{27}\log a+\tfrac49\left(\tan^{-1}\left(\sqrt{\a}\right)\right)^2
-\frac{2\{3+\a(3-4\a)\}\tan^{-1}\left(\sqrt{\a}\right)}{9\a^{5/2}}+\frac{6+4\a}{9\a^2}-\frac{44}{45}\,.
\end{align*}
\eop

\section{Appendix}
\label{sec:appendix}
\setcounter{equation}{0}
{\bf Proof of (\ref{relevant_exp}):}
By (4.6) in \cite{cabogr:94} we have, as a consequence of Lemma 2.4 on page 339 of
\cite{gr:88},
\begin{align*}
EU(1)^4V(s)^4=&\tfrac45EU(1)^4\left\{V(1)^4\exp\left\{-\tfrac12\s V(1)^2\right\}
-\frac{V(1)^2}{\s}\exp\left\{-\tfrac12\s V(1)^2\right\}\right.\\
&\qquad\qquad\qquad\qquad\qquad\qquad\left.+\frac2{\s^2}\left\{1-\exp\left\{-\tfrac12\s V(1)^2\right\}\right\}\right\},
\end{align*}
where $\s=s-1$. Again using Lemma 2.4 in \cite{gr:88} we have:
$$
EU(1)^4\exp\left\{-\tfrac12\s V(1)^2\right\}
=\int_{\R_+^2}x^4\exp\left\{-\tfrac12\s y^2-\tfrac12(x+y)^2\right\}\,dxdy.
$$
Let the integral $\int_{\R_+^2}x^4\exp\left\{\tfrac12\s y^2-\tfrac12(x+y)^2\right\}\,dxdy$ be denoted by $\f(\s)$. Then
$$
\f(\s)=\frac{3(1+\s)^2\tan^{-1}\left(\sqrt{\s}\right)}{\s^{5/2}}-\frac{3+5\s}{\s^2}\,.
$$
This can be verified by using the recursive relations of Lemma 4.1 in combination with (4.8) in \cite{cabogr:94}. Note that
$$
\f(0)\stackrel{\mbox{\small def}}=\lim_{\s\downarrow0}\f(\s)=\frac85\,.
$$
Furthermore,
$$
EU(1)^4V(1)^2\exp\left\{-\tfrac12\s V(1)^2\right\}=-2\f'(\s)=
\frac{3(1+\s)(5+\s)\tan^{-1}\left(\sqrt{\s}\right)}{\s^{7/2}}-\frac{15+13\s}{\s^3}\,,
$$
and
$$
EU(1)^4V(1)^4\exp\left\{-\tfrac12\s V(1)^2\right\}=4\f''(\s)=\frac{3\{35+3\s(10+\s)\}
\tan^{-1}\left(\sqrt{\s}\right)}{\s^{9/2}}-\frac{5\s(21+11\s)}{\s^4}\,.
$$
This yields:
\begin{equation}
EU(1)^4V(s)^4\longrightarrow\frac{64}{105}\,,\s\downarrow0,
\end{equation}
which coincides with (4.9) in \cite{cabogr:94}.

For the other term we have:
\begin{align*}
&EU(1)^3V(1)V(s)^4\\
&=\frac45E U(1)^3\left\{V(1)^5\exp\left\{-\tfrac12\s V(1)^2\right\}
-\frac{V(1)^3\exp\left\{-\tfrac12\s V(1)^2\right\}}{\s}\right.\\
&\qquad\qquad\qquad\qquad\qquad\qquad\qquad\qquad\qquad\qquad+\left.\frac{2V(1)\left(1-\exp\left\{-\tfrac12\s
V(1)^2\right\}\right)}{\s^2}\right\},
\end{align*}
where $\s=s-1$. We have:
$$
EU(1)^3V(1)\exp\left\{-\tfrac12\s V(1)^2\right\}
=\frac2{\s}+\frac3{\s^2}-\frac{3(1+\s)\tan^{-1}\left(\sqrt{\s}\right)}{\s^{5/2}}
$$
Define $\psi(s)=EU(1)^3V(1)\exp\left\{-\tfrac12\s V(1)^2\right\}$. Then
$$
EU(1)^3V(1)^3\exp\left\{-\tfrac12\s V(1)^2\right\}=-2\psi'(\s)=
\frac{15+4\s}{\s^3}-\frac{3(5+3\s)\tan^{-1}\left(\sqrt{\s}\right)}{\s^{7/2}},
$$
and
$$
EU(1)^3V(1)^5\exp\left\{-\tfrac12\s V(1)^2\right\}=4\psi''(\s)=
\frac{105+\s(115+16\s)}{\s^4(1+\s)}-\frac{15(7+3\s)\tan^{-1}\left(\sqrt{\s}\right)}{\s^{9/2}}.
$$
Hence
\begin{align*}
&EU(1)^3V(1)V(s)^4\\
&=\frac{8\left\{210+\s(215+2\s(11+\s))\right\}}{25\s^4(1+\s)}-\frac{24(14+5\s)\tan^{-1}\left(\sqrt{\s}\right)}{5\s^{9/2}}
\longrightarrow \frac{16}{21},\,\s\downarrow0.
\end{align*}
Formula (\ref{relevant_exp}) now follows.\eop

\section{Concluding remarks}
\label{sec:conclusion}
There is a remarkable analogy between the behavior of the left-lower convex hull of the Poisson point process, discussed above, and the least concave majorant of (one-sided) Brownian motion without drift, as analyzed in \cite{piet:83}. In the same way there is an analogy between the behavior of the lower convex hull of the Poisson point process inside a parabola, as analyzed in \cite{gr:88} and \cite{nagaev:95}, and the least concave majorant of Brownian motion with a parabolic drift, as analyzed in \cite{piet:89} and \cite{piet:11b}. Why this is the case is still somewhat of a mystery and deserves (in my view) further investigation.

\vspace{0.4cm}
\noindent
{\bf Acknowledgements} I want to thank Tomasz Schreiber for sending me the unpublished preprint \cite{nagaev_kham:91} and Christian Buchta for pointing out some typographical errors in the original version of this manuscript (written in 2006).

\end{document}